\documentclass{llncs}

\usepackage[utf8]{inputenc}
\usepackage[T1]{fontenc}
\usepackage[brazilian]{babel}
\usepackage{csquotes}
\usepackage{hyphenat}
\usepackage{amssymb}
\usepackage{amsmath}
\usepackage{array}

\usepackage[locale=FR]{siunitx}

\usepackage{graphicx}
\graphicspath{ {images/} }

\usepackage{listings}

\begin{document}

\title{On Coupled Dirac Systems under Boundary Condition}
\titlerunning{Experimentos com Algoritmos Genéticos}

\author{ Xu Yang* and Xin Li}
%
%
%
\institute{School of Mathematics , Yunnan Normal University\\
\email{yangxu02465129@163.com}}

\maketitle              
\renewcommand\refname{References}
\renewcommand{\abstractname}{Abstract}
\begin{abstract}
In this article we study the existence of solutions for the Dirac systems
\begin{equation}\label{e:0.1}
\left\{
\begin{array}{c}
  Pu=\frac{\partial H}{\partial v}(x,u,v) \quad\hbox{on} \ M,\\
  Pv=\frac{\partial H}{\partial u}(x,u,v) \quad\hbox{on} \ M,\\
  B_{\text{CHI}}u= B_{\text{CHI}}v=0\quad\hbox{on} \ \partial M
\end{array}
\right.
\end{equation}
where $M$ is an $m$-dimensional compact oriented Riemannian spin manifold with smooth boundary $\partial M$, $P$ is the Dirac operator under the boundary condition $B_{\text{CHI}}u= B_{\text{CHI}}v=0$ on $\partial M$, $ u,v\in C^{\infty}(M,\Sigma M)$ are spinors. Using an analytic framework of proper products of fractional  Sobolev spaces, the solutions existence results of the coupled Dirac systems are obtained for nonlinearity with superquadratic growth rates.
\end{abstract}
\section{Introduction and main results}
Dirac operators on compact spin manifolds play prominent role in the geometry and mathematical physics, such as the generalized Weierstrass representation of the surface in three manifolds \cite{Fri1} and the supersymmetric nonlinear sigma model in quantum field theory \cite{CJLW1,CJLW2}. The existence of solutions of the Dirac equation has been studied on compact spin manifolds without boundaries by Ammann\cite{Am1,Am2}, Raulot\cite{Rau},Isobe\cite{Iso1,Iso2,Iso3},Gong and Lu \cite{GoL1,GoL2},Yang\cite{YJL}. In addition, be different with these existing works, Ding and Li \cite{DiL} studied  a class of boundary value problem on a compact spin manifold $M$ with smooth boundary. The problem is a general relativistic model of confined particles by means of nonlinear Dirac fields on $M$. In this paper, we are concerned with a nonlinear Dirac systems on compact spin manifolds with smooth boundary, and deal with some existence results.
\par
Spin manifold $(M,g)$  equipped with a spin structure
$\sigma: P_{spin(M)}\rightarrow P_{so(M)}$, and let $\Sigma M= P_{spin(M)}\times_{\sigma}\Sigma_{m}$ denote the complex spinor bundle on $M$, which is a complex vector bundle of rank $2^{[m/2]}$ endowed with the spinorial Levi-Civita connection $\nabla$ and a pointwise Hermitian scalar product $\langle\cdot,\cdot\rangle$. We always assume $m\ge 2$ in this paper. Consider Whitney direct sum $\Sigma M\oplus\Sigma M$ and Write a point of it as $(x,u,v)$, where $x\in M$ and $u,v\in \Sigma_x M$. $P$ is the Dirac operator under the boundary condition $B_{\text{CHI}}u= B_{\text{CHI}}v=0$ on $\partial M$. We consider the following system of the coupled equations
\begin{equation}\label{e:1.1}
\left\{
\begin{array}{c}
  Pu=\frac{\partial H}{\partial v}(x,u,v) \quad\hbox{on} \ M,\\
  Pv=\frac{\partial H}{\partial u}(x,u,v) \quad\hbox{on} \ M,\\
  B_{\text{CHI}}u= B_{\text{CHI}}v=0\quad\hbox{on} \ \partial M
\end{array}
\right.
\end{equation}
where  fiber preserving map $H:\Sigma M\oplus  \Sigma M \rightarrow \mathbb{R}$ is a real valued superquadratic function of class $C^{1}$ with subcritical growth rates.
(\ref{e:1.1}) is the Euler-lagrange equation of the functional
\begin{equation}\label{e:1.2}
 \mathfrak{L}(u,v)=\int_{M}\{\langle Pu,v\rangle-H(x,u,v)\}dx,
\end{equation}
where $dx$ is the Riemann volume measure on $M$ with respect to the metric $g$, $\langle\cdot,\cdot\rangle$ is the compatible metric on $\Sigma M$.
\par
Let $z=(u,v)$,
$L=\left(
  \begin{array}{cc}
    0 & P \\
    P & 0 \\
  \end{array}
\right)$, and  $Lz:=(Pv,Pu)$. Then (\ref{e:1.2}) becomes
\begin{equation}\label{e:1.3}
 \mathfrak{L}(z)=\int_{M}\left\{\frac{1}{2}\langle Lz,z\rangle-H(x,z)\right\}dx.
\end{equation}

The problem (\ref{e:1.1}) describes two coupled fermionic fields in quantum. It can be viewed as a spinorial analogue of other strongly indefinite variational problems such as infinite dynamical systems \cite{BaDi1,BaDi2} and
elliptic systems \cite{AnVo,FF}. A typical way to deal with such problems is the min-max method of Benci and Rabinowitz \cite{BeRa}, including the mountain pass theorem, linking arguments. Another is a homological method,  Morse theory and  Rabinowitz-Floer homology as in \cite{Abb,BaLi,KrSz,Maa,MaM}.
 For the Dirac operator associated with appropriate boundary condition , Farinell, Schwarz \cite{FaSc} prove that Dirac operator $P$ is elliptic and extends to a self-adjoint operator with  a discrete spectrum. Except for the work, this type of boundary condition was also considered  by Hijazi, Montiel,Rold$\acute{a}$n, Zhang\cite{HMR,HMZ}. In this paper we use the techniques introduced by Hulshof and Van Der Vorst \cite{HV} to prove the existence of solutions of (\ref{e:1.1}), and apply a generalized fountain theorem established by Batkam and Colin \cite{BaCo} to obtain infinitely many solutions of the coupled Dirac system provided the nonlinearity $H$ is even.
\par
In the following we assume that two real numbers $p,q$ satisfy
$$
\frac{1}{p}+\frac{1}{q}>\frac{m-1}{m}.
$$
For the nonlinearity $H$, we make the following hypotheses:
\par
$H\in C^0(\Sigma M\oplus \Sigma M, \mathbb{R})$ is $C^1$ in the fiber direction.
 Real constants $2<p,q< 2^{\ast}:=\frac{2m}{m-1}$ and $\mu$ satisfy
\begin{description}
\item[(i)] $\max \{p-1,q-1,2\}<\mu\leq\min \{p,q\}$,
\item[(ii)] $(1-\frac{1}{p})\max\{\frac{p}{\mu},\frac{q}{\mu}\} <\frac{1}{2}+\frac{1}{2m}$,
\item[(iii)] $(1-\frac{1}{q})\max\{\frac{p}{\mu},\frac{q}{\mu}\} <\frac{1}{2}+\frac{1}{2m}$,
\item[(iv)] $\frac{p-1}{p}\frac{q}{\mu}<1$ and $\frac{q-1}{q}\frac{p}{\mu}<1$.
\end{description}
\par
\noindent{(${\bf H}_1$)} There exist a constant $C_1>0$ such that
\begin{eqnarray}
  \left| H_{u}(x,u,v)\right|&\leq& C_1\left(|u|^{p-1}+|v|^{\frac{(p-1)q}{p}}+1\right),\label{e:1.4} \\
  \left| H_{v}(x,u,v)\right|&\leq& C_1 \left(|v|^{q-1}+|u|^{\frac{(q-1)p}{q}}+1\right).\label{e:1.5}
\end{eqnarray}
\\
\noindent{(${\bf H}_2$)} There exist $R_{1}>0$, such that
\begin{equation}\label{e:1.6}
 0<\mu H(x,u,v)\leq \langle H_u,u\rangle+\langle H_v,v\rangle,
\end{equation}
for all $(x,u,v)\in \Sigma M\oplus\Sigma M$ with $|(u,v)|\geq R_1$.
\\
\noindent{(${\bf H}_3$)} $H(x,u,v)\geq 0$ for all $(x,u,v)\in \Sigma M\oplus\Sigma M$.
\\
\noindent{(${\bf H}_4$)} $H(x,u,v)=o(|(u,v)|^2)$ as $|(u,v)|\rightarrow 0$ uniformly for $x\in M$.
\\
\noindent{(${\bf H}_5$)} $H(x,-u,-v)=H(x,u,v)$ for any $(x,u,v)\in \Sigma M\oplus\Sigma M$.
\par
{\bf Notice:} That $H(x,u,v)=|u|^p+|v|^q$ satisfies these conditions. In \cite{Yan} existence results for the Dirac system without boundary condition are give under the same assumptions on $H(x,u,v)$. The above condition (${\bf H}_2$) looses associated condition in  Gong and Lu \cite{GoL1}.
\par
Our main result is as follow.
\begin{theorem}\label{the:1.1}
If the above $H$ satisfies $({\bf H}_1)-({\bf H}_4)$, then  Dirac system(\ref{e:1.1}) possesses at least one solution.
\end{theorem}

Furthermore, for odd nonlinearities we have the following multiplicity result:
\begin{theorem}\label{the:1.2}
If the above $H$ satisfies $({\bf H}_1)-({\bf H}_5)$,
then there exists a sequence of solutions $\{(u_k,v_k)\}_{k=1}^\infty $ to (\ref{e:1.1}) with $\mathfrak{L}(u_k,v_k)\rightarrow \infty$ as $k\rightarrow \infty$.
\end{theorem}

\section{About boundary condition}
We collect here some basic definitions and facts about spin structures on manifolds and Dirac operators. For more detailed exposition, please consult\cite{DiL,LaM,Fri2}.
\par
 Define $(\Sigma M,\langle ,\rangle,\gamma,\nabla)$ is a Dirac bundle if $\nabla :C^\infty
 (M,\Sigma M)\rightarrow C^\infty(M,T^\ast M \otimes\Sigma M)$ and $\gamma: C^\infty(M,T^\ast M \otimes\Sigma M) \rightarrow C^\infty
 (M,\Sigma M) $ satisfies:
 $$\nabla \gamma=0,\langle \gamma (\omega)\phi,\psi\rangle +\langle\phi, \gamma (\omega)\psi\rangle=0,\nabla(\gamma (\omega)\phi)=
 \gamma (\nabla\omega)\phi+\gamma (\omega)\nabla\phi$$
for any $\omega\in TM$ and $\phi,\psi\in C^\infty(M,\Sigma M)$.
We have used the identification $T^\ast M\cong TM $ by the metric on $M$, then $C^\infty(M,T^\ast M \otimes\Sigma M)\cong C^\infty(M,T M \otimes\Sigma M)$,
 therefore, Dirac operator $P$ act on spinors on $M$ is defined  by $$P\psi=\gamma\circ\nabla \psi, \quad for\ any\ \psi \in C^\infty(M,\Sigma M).$$
In particular, if we choose an local orthogonal tangent frame $\{e_1,e_2,...,e_m\}$, the Dirac operator $P$ becomes
$$
P\psi=\Sigma_{j=1}^m \gamma(e_j)\nabla_{e_j}\psi, \quad for\ any\ \psi \in C^\infty(M,\Sigma M).
$$
Then we  consider a Chirality operator associated with the Dirac bundle
 $$
 \langle \Sigma M,\langle\cdot\rangle,\gamma,\nabla\rangle.
 $$
 If a linear map
 $F:End_\mathbb{C}(\Sigma M)\rightarrow End_\mathbb{C}(\Sigma M)$ satisfies
 $$
 F^2=\text{Id},\langle F\psi,F\varphi\rangle=\langle \psi,\varphi\rangle,\nabla_X(F\psi)=F\nabla_X\psi,\gamma(X)F\psi=-F\gamma(X)\psi
 $$
for each vector field $X\in TM$ and spinor fields $\psi,\varphi\in C^\infty(M,\Sigma M)$.
\par
The boundary hypersurface $\partial M$ is also a spin manifold and so we have the corresponding spinor bundle $\Sigma \partial M$, the clifford multiplication  $\gamma^{\partial M}$, the spin connection $\nabla^{\partial M}$ and the intrinsic Dirac operator $P^{\partial M}$. In \cite{HMR},  Hijazi,  Montiel, Rold$\acute{a}$n show that  the restricted Hermitian bundle $\Sigma M\mid_{\partial M}$ can be  identified with the intrinsic Hermitian spinor bundle $\Sigma \partial M$, provided that $m$ is odd. Instead, if $m$ is even, the restricted Hermitian bundle $\Sigma M\mid_{\partial M}$ could be identified with the sum  $\Sigma \partial M\oplus \Sigma \partial M$.

Define an operator $\Gamma:=F\mid_{\partial M}\gamma(N)$, where $N$ denotes the unite inner normal vector field on $\partial M$.
By the definition , we know $F$ is a local operator on the spinor bundle over  $\partial M$. $F$ is a self-adjoint operator and has two eigenvalues $+1$ and $-1$. The corresponding eigenspaces are
$$
\Gamma_+=\{\phi\in C^\infty(\partial M,\Sigma\partial M)| F\phi=\phi\},\Gamma_-=\{\phi\in C^\infty(\partial M,\Sigma\partial M)| F\phi=-\phi\}.
$$
Now consider the boundary condition in problem (\ref{e:1.1}).The operator $B_{\text{CHI}}$ on $\partial M$ is defined as $B_{\text{CHI}}=\frac{1}{2}(\text{Id}-F)=\frac{1}{2}(\text{Id}-\gamma(N)G)$. Then $B_{\text{CHI}}$ is a self-adjoint operator and it is not difficult to see that
$$
B_{\text{CHI}}\psi\mid_{\partial M}=0\Leftrightarrow \psi\in  C^\infty( M,\Sigma M) \quad and  \quad\psi\mid_{\partial M}\in \Gamma_+.
$$
For the space $C^\infty( M,\Sigma M)$, define an inner product
$$
(\psi,\varphi)_{H^1}=(\psi,\varphi)_2+(\nabla\psi,\nabla\varphi)_2, \quad \forall\psi,\varphi\in C^\infty( M,\Sigma M)
$$
Then $H^1( M,\Sigma M)$ is the completion of the space $ C^\infty( M,\Sigma M)$ with respect to the norm $\parallel\cdot\parallel_{H^1}$.
Since $P$ is a first operator, it extends to a linear operator $P:H^1( M,\Sigma M)\rightarrow L^2( M,\Sigma M)$ and $P|_{\partial M}:H^1( M,\Sigma M)\rightarrow L^2(\partial M,\Sigma\partial M)$.Let
$$
\mathfrak{D}(P)=\{\psi\in H^1( M,\Sigma M)\mid\psi\mid_{\partial M}\in L^2\Gamma_+\}.
$$
Then the Dirac operator $P$ with Chirality boundary condition $B_{\text{CHI}}\psi\mid_{\partial M}=0$ is well defined in the domain $\mathfrak{D}(P)$.
For simplicity, in the following, we will denote the $\mathfrak{D}(P)$ by $\mathfrak{D}$.
\par
For $\psi,\varphi\in \mathfrak{D}$,By the integrated version of Lichnerowitz Formula, we have $(P\psi,\varphi)=(\psi,P\varphi)$,Actually, $P$ is a self-adjoint operator in $L^2(M,\Sigma M)$ with domain $\mathfrak{D}$.

\section{The analytic framework}
Recall the operator $P$ is an unbounded self-adjoint operator, Let \text{Spec}(P) denote the spectrum of the Dirac operator $P$ under the local boundary condition  $B_{\text{CHI}}\psi\mid_{\partial M}=0$.
Moreover, there exists a complete orthonormal basis $\{\eta_{k}\}_{k=1}^{\infty}$ of the Hilbert space $\mathfrak{D}$  consisting of the eigenspinors of the operator $P: P\eta_{k}=\lambda_{k}^{\text{CHI}}\eta_{k}$. Moreover, $|\lambda_{k}^{\text{CHI}}|\rightarrow \infty$ (as $k\rightarrow \infty$) and all corresponding eigenvalues $\{\lambda_{k}^{\text{CHI}}\}_{k=1}^{\infty}$ have finite multiplicity.
\par
If $(M,g)$ has positive scalar curvature, it is obviously $0\notin \text{Spec}(P)$ By the Fridrich's inequality.
\par
Define
$$
(\psi,\varphi)_{1,2}:=(P\psi,P\varphi)_2,\quad \forall \psi,\varphi \in \mathfrak{D}.
$$
where $(\cdot ,\cdot )_2$ is the $L^2$-inner product on spinors.
We denote $\parallel \psi\parallel_{1,2}=(\psi,\psi)_{1,2}^\frac{1}{2}$.
Since $P$ is self-adjoint, we can prove $ (\mathfrak{D},\parallel \cdot\parallel_{1,2})$ is a Hilbert space. In addition, the following lemma tells us that
$ (\mathfrak{D},\parallel \cdot\parallel_{1,2})$ embeds continuously in $H^1( M,\Sigma M)$.
\par
\begin{lemma}\cite{FaSc}\label{the:3.1}
There exists a universal constant $C\in [1,+\infty]$ such that
$$
\parallel\psi\parallel^2_{H^1}\leq C \parallel P\psi\parallel^2_2,\quad \psi\in \mathfrak{D}.
$$
\end{lemma}
Let $\mid P\mid$ denote the absolute value operator of $P$ defined in $L^2(M,\Sigma M)$, and $\mid P\mid^{\frac{1}{2}}$ is the square root operator of $\mid P\mid$ with its spectrum
$$
\text{Spec}(\mid P\mid^{\frac{1}{2}})=\{\mid \lambda_k^{\text{CHI}}\mid^{\frac{1}{2}} \ \mid \ k=1,2,3,\cdot\cdot\cdot\}.
$$
The inner product is define by
$$
(\psi,\varphi)_{\frac{1}{2},2}:=(|P|^{\frac{1}{2}}\psi,|P|^{\frac{1}{2}}\varphi)_2,\quad \forall \psi,\varphi \in \mathfrak{D}.
$$
\par
Let $E^{\frac{1}{2}}:=\mathfrak{D}(|P|^\frac{1}{2})$ denote the domain of the operator $\mid P\mid^{\frac{1}{2}}$ and $\mid P\mid^{\frac{1}{2}}$ is self-adjoint in $L^2(M,\Sigma M)$. We denote $\parallel \psi\parallel_{\frac{1}{2},2}=(\psi,\psi)_{\frac{1}{2},2}^\frac{1}{2}$ . Then  $(E^{\frac{1}{2}},\parallel \cdot \parallel_{\frac{1}{2},2})$ is a Hilbert space, Its dual space is denoted by  $E^{-\frac{1}{2}}$. Then $\mid P\mid^{-1}$ is a Hilbert space isomorphism from  $E^{-\frac{1}{2}}$ to $E^{\frac{1}{2}}$ with respect to the equivalent inner products  $(\cdot,\cdot)_{\frac{1}{2},2}$.
\par
Consider the Hilbert space
$$E:=E^\frac{1}{2}\times E^\frac{1}{2}$$
with inner product
$$
((u_1,v_1),(u_2,v_2))_E:=(u_1,u_2)_{\frac{1}{2},2}+( v_1,v_2)_{\frac{1}{2},2}\quad \forall (u_1,v_1),(u_2,v_2)\in E
$$
and norm $ \parallel z\parallel_E:=(\parallel u\parallel_{\frac{1}{2},2}^2+\parallel v\parallel_{\frac{1}{2},2}^2)^\frac{1}{2}$ for $z=(u,v)\in E$.
Let $E^\ast=E^{-\frac{1}{2}}\times E^{-\frac{1}{2}}$, which is the dual space of $E$. Then
$$
\left(
  \begin{array}{cc}
    \mid P\mid^{-1} & 0 \\
    0 & \mid P\mid^{-1} \\
  \end{array}
\right):E^\ast\rightarrow E
$$
is a Hilbert space isomorphism by the arguments above. Hence

$$
\left(
  \begin{array}{cc}
    \mid P\mid^{-1} & 0 \\
    0 & \mid P\mid^{-1} \\
  \end{array}
\right)\circ
\left(
  \begin{array}{cc}
    0 & P \\
    P & 0 \\
  \end{array}
\right)
=
\left(
  \begin{array}{cc}
    0& \mid P\mid^{-1}P\\
    \mid P\mid^{-1}P & 0 \\
  \end{array}
\right)
,
$$
denote
$$
B=\left(
  \begin{array}{cc}
    0& \mid P\mid^{-1}P\\
    \mid P\mid^{-1}P & 0 \\
  \end{array}
\right).
$$
\par
It is a self-adjoint isometry operator and $B\circ B=\text{Id}: E\rightarrow E$ is identity operator. Introducing the "diagonals"
$$
E_\pm=\{(\pm \mid P\mid^{-1}P v,v)\mid v\in E^\frac{1}{2}\},
$$
We have
$$
E=E_+\oplus E_-=\{z=z^++z^- \mid z^\pm\in E_\pm\}.
$$
Note that $B z^\pm=\pm z^\pm$, so that $E_+$ and $E_-$ are  the mutually orthogonal eigenspaces of the eigenvalues $1$ and $-1$ of $B$. Orthonormal bases consisting of eigenvectors of $E_\pm$ are given by
$$
\{\frac{1}{\sqrt 2}(|\lambda_k^{\text{CHI}}|^{-\frac{1}{2}}\eta_k,\pm|\lambda_k^{\text{CHI}}|^{-\frac{3}{2}}\lambda_k^{\text{CHI}}\eta_k)\}_{k=1}^{\infty}.
$$
Then
$$
 E_+:= \overline{\bigoplus\limits_{j=1}^{\infty}\mathbb{R}e_j} \quad with\  e_j=\frac{1}{\sqrt 2}(|\lambda_k^{\text{CHI}}|^{-\frac{1}{2}}\eta_k,|\lambda_k^{\text{CHI}}|^{-\frac{3}{2}}\lambda_k^{\text{CHI}}\eta_k).
$$
We have
\begin{eqnarray}
(L z,z)_2&=&\int_{M}\langle L z,z\rangle dx=2\int_{M}\langle Pu,v\rangle dx=2\int_{M}\langle \mid P\mid^{\frac{1}{2}}\mid P\mid^{-1}Pu,\mid P\mid^{\frac{1}{2}}v\rangle dx \nonumber\\
&=&2(\mid P\mid^{-1}Pu,v)_{\frac{1}{2},2}=(Bz,z)_E.\nonumber
\end{eqnarray}
Then for each $z=z^++z^-$, we have
$$
(L z,z)_2=\int_{M}\langle L z,z\rangle dx=(Bz,z)_E=\parallel z^+\parallel^2_E-\parallel z^-\parallel^2_E.
$$
\par
Now we can define a functional $\mathcal{L}:E\rightarrow \mathbb{R}$ as
\begin{equation}\label{e:3.1}
 \mathfrak{L}(z)=\frac{1}{2}(Bz,z)_{\frac{1}{2},2}+\mathcal{H}(z)=\frac{1}{2}\parallel z^+\parallel^2_E-\frac{1}{2}\parallel z^-\parallel^2_E-\mathcal{H}(z),
\end{equation}
\par
 where $\mathcal{H}(z)=\int_M H(x,z)dx$.
 \par
 Since $M$ is compact, by the assumption $(H_1)$ and integrating we obtain
\begin{equation}\label{e:3.2}
  \mid H(x,u,v)\mid\leq \mid H(x,0,v)\mid +C_1\left(|u|^{p-1}+|v|^{\frac{(p-1)q}{p}}+1\right)\mid u\mid
\end{equation}
and let $u=0$, similarly, we prove that
\begin{equation}\label{e:3.3}
  \mid H(x,0,v)\mid\leq \mid H(x,0,0)\mid+C_1\left(\mid v\mid^{q}+\mid v\mid +1\right).
\end{equation}
From $(\ref{e:3.2}),(\ref{e:3.3})$ and Young' inequality to derive
\begin{equation}\label{e:3.4}
  \mid H(x,u,v)\mid\leq C\left(|u|^{p}+|v|^{q}+1\right),
\end{equation}
for some constant $C>0$.
\par
By an analysis of interpolation of the Sobolev spaces,
\par
$$E^\frac{1}{2}=[\mathfrak{D},L^2]_{\frac{1}{2}}\quad \text{and} \quad H^\frac{1}{2}(M,\Sigma M)=[H^1,L^2]_\frac{1}{2},$$
where $[\cdot,\cdot]_\frac{1}{2}$ means the interpolation and $H^\frac{1}{2}(M,\Sigma M)$ coincides with the usual $L^2-$ Sobolev space of order $s$, $W^{s,2}(M,\Sigma M)$,see \cite{Ada}.
 \par
 Since $ (\mathfrak{D},\parallel \cdot\parallel_{1,2})$ embeds $H^1( M,\Sigma M)$ continuously, there holds the continuous embedding
$$
E^\frac{1}{2}\times E^\frac{1}{2}=[\mathfrak{D},L^2]_{\frac{1}{2}}\times [\mathfrak{D},L^2]_{\frac{1}{2}}\hookrightarrow [H^1,L^2]_{\frac{1}{2}}\times [H^1,L^2]_{\frac{1}{2}}=H^{\frac{1}{2}}\times H^{\frac{1}{2}}.
$$
\par
By the fact $H^\frac{1}{2}(M,\Sigma M)$ embeds $L^p(M,\Sigma M)$
continuously for $1\leq p\leq\frac{2m}{m-1}$. Moreover, this embedding is compact if $1\leq p<\frac{2m}{m-1}$.
we obtain the following lemma:
\begin{lemma}\label{the:3.2}
$E^\frac{1}{2}\times E^\frac{1}{2}$ embeds  $L^p( M,\Sigma M)\times L^q( M,\Sigma M)$ continuously for $1\leq p\leq \frac{2m}{m-1},1\leq q\leq \frac{2m}{m-1}$. Moreover, this embedding is compact if $1\leq p< \frac{2m}{m-1},1\leq q< \frac{2m}{m-1}$. Where norm $\parallel z\parallel_{L^p\times L^q}:=(\parallel u\parallel_p^p+\parallel v\parallel_q^q)^\frac{1}{2}$ for $z=(u,v)\in L^p( M,\Sigma M)\times L^q( M,\Sigma M)$.
\end{lemma}

Then  using $(\ref{e:3.4})$ we can define the functional $\mathcal{H}:E\rightarrow \mathbb{R}$ as
$$
\mathcal{H}(u,v)=\int_M H(x,u,v)dx,
$$
is of class $C^1$ and its derivative at $(u,v)\in E$ is given by
$$
D\mathcal{H}(u,v)(\xi,\zeta)=\int_M( \langle H_u(x,u,v),\xi\rangle +\langle H_v(x,u,v),\zeta\rangle) dx \quad (\xi,\zeta)\in E.
$$
Moreover $D \mathcal{H}:E\rightarrow E^*$ is a compact operator.
\par
In fact, Using H$\ddot{o}$lder inequality and embeddings we have
$$
\int_M \langle H_u(x,u,v),\xi\rangle  dx\leq C(\parallel u\parallel_{\frac{1}{2},2}^{p-1}+\parallel v\parallel_{\frac{1}{2},2}^\frac{(p-1)q}{p}+1)\parallel\xi\parallel_{\frac{1}{2},2}.
$$
In a similar way we obtain an inequality for the derivative with respect to $v$. Thus $D\mathcal{H}(u,v)$ is well defined and bounded in $E$.
Next, by the Sobolev embeddings, usual arguments give that  $D\mathcal{H}(u,v)$ is compact.

\section{The Palais-Smale condition for $\mathfrak{L}$}
Let $F$ be a $C^1$ functional  on a Banach space $E$, $c\in \mathbb{R}$. Recall that
 a sequence $\{x_{n}\}\subset E$ is called a $(PS)_{c}-$ sequence if
 $F(x_{n})\rightarrow c$ as $n\rightarrow \infty$ and
 $\parallel DF(x_{n})\parallel_{E^{*}}\rightarrow 0$ as $n\rightarrow \infty$.
If all $(PS)_{c}$ -sequences converge in $E$, we say that $F$ satisfies the $(PS)_{c}$ condition.
In this section we prove the$(PS)_{c}$ condition for $\mathfrak{L}$.
\par

\begin{lemma}\label{lem:4.1}
Suppose $H$ satisfies $({\bf H}_1)$,$({\bf H}_2)$. Then for any $c\in \mathbb{R}$, $\mathfrak{L}$ satisfies the $(PS)_{c}$-condition with respect to $E$.
\end{lemma}
\noindent{\bf Proof}.\quad  let $\{z_n\}=\{(u_n,v_n)\}\subset E$ be a $(PS)_{c}$-sequence with respect to $E$, i.e., $z_n\in E$ and satisfy
\begin{equation}\label{e:4.1}
  \mathfrak{L}(z_n)\rightarrow c \quad \hbox{as} \ n\rightarrow \infty
\end{equation}
and
\begin{equation}\label{e:4.2}
  \parallel D\mathfrak{L}(z_n)\parallel_{E^*}\rightarrow 0\quad
   \hbox{as} \ n \rightarrow\infty.
\end{equation}
\noindent{\it Claim }.  $\{z_n\}\subset E$ is bounded.

The condition (${\bf H}_2$) implies that there are constants $C_2, C_3>0$ such that
\begin{equation}\label{e:4.3}
  H(x,u_n,v_n)\geq C_2(|u_n|^\mu+|v_n|^{\mu})-C_3.
\end{equation}
 See \cite{Fel} for a proof.
By (\ref{e:4.1})-(\ref{e:4.3}) and (${\bf H}_2$) ,  for large $n$ we have
\begin{eqnarray}\label{e:4.4}
C+\parallel z_n\parallel_{E}&\geq& 2\mathfrak{L}(z_n)-\langle D\mathfrak{L}(z_n),z_n\rangle\nonumber\\
&=& \int_M\langle H_z,z_n\rangle dx-2\int_M H(x,z_n)dx\nonumber\\
&\geq&(\mu-2)\int_{M}H(x,z_n)dx-C\nonumber\\
&\geq& C_2(\mu-2)\int_{M}(|u_n|^\mu+|v_n|^{\mu})dx-C_3(\mu-2)\mid M\mid-C.
\end{eqnarray}
Hereafter  $C$ denote various positive constants which do not depend on $n$. Clearly, (\ref{e:4.4}) implies that for large $n$,
$\|u_n\|_\mu^\mu+\|v_n\|_\mu^\mu\leq C(1+\|z_n\|_E)$ and so
\begin{equation}\label{e:4.5}
  \|u_n\|_\mu^\mu\leq C(1+\|z_n\|_E)\quad\hbox{and}\quad
  \|v_n\|_\mu^\mu\leq C(1+\|z_n\|_E).
  \end{equation}

Write $z_n=z_n^{-}+z_n^{+}$ according to the decomposition $E=E_{-}\oplus  E_{+}$.
By (\ref{e:4.2}), for large $n$, we have
\begin{eqnarray}\label{e:4.7}
  \left| \langle D\mathfrak{L}(z_n),z_n^+\rangle\right| &=&\left|\|z_n^+\|^2_E-\int_M\langle z_n^+,H_z\rangle dx\right|
  \leq \|z_n^+\|_E.
\end{eqnarray}
This  and (${\bf H}_1$) lead to
\begin{eqnarray}\label{e:4.8}
&&\|z_n^+\|^2_E \leq\left| \int_M\langle z_n^+,H_z\rangle dx \right|+\|z_n^+\|_E\nonumber \\
  &\le&\left| \int_M\langle u_n^+,H_u\rangle dx \right|+\left| \int_M\langle v_n^+,H_v\rangle dx \right|+\|z_n^+\|_E\nonumber \\
 &\leq&  \int_M| u_n^+||H_u| dx + \int_M | v_n^+||H_v| dx+\|z_n^+\|_E\nonumber\\
 &\leq& \int_M c_1\left(|u_n|^{p-1}+|v_n|^{\frac{(p-1)q}{p}}+1\right)| u_n^+|dx+\nonumber\\
 &+&\int_M c_1 \left(|v_n|^{q-1}+|u_n|^{\frac{(q-1)p}{q}}+1\right)| v_n^+|+\|z_n^+\|_E.
\end{eqnarray}
 Note that
 $1<\frac{\mu}{\mu-p+1}<2^\ast$ and $1<\frac{\mu p}{\mu p-(p-1)q}<2^\ast$
 by the conditions (ii) and (iv) above (${\bf H}_1$). We derive
 $$
 E^\frac{1}{2}\hookrightarrow L^\frac{\mu}{\mu-p+1},\quad E^\frac{1}{2}\hookrightarrow L^\frac{\mu p}{\mu p-(p-1)q},
 $$
 and therefore
 \begin{eqnarray}
   \int_M |u_n|^{p-1}|u_n^+|dx &\leq& C \| u_n\|_{\mu}^{p-1}\|u_n^+\|_{\frac{1}{2},2},\label{e:4.9}\\
    \int_M |v_n|^\frac{(p-1)q}{p}|u_n^+|dx&\leq& C \| v_n\|_{\mu}^\frac{(p-1)q}{p}\|u_n^+\|_{\frac{1}{2},2}.\label{e:4.10}
 \end{eqnarray}
Using the conditions (iii) and (iv) above (${\bf H}_1$),  an analogous reasoning yields
 \begin{eqnarray}
   \int_M |v_n|^{q-1}|v_n^+|dx &\leq& C \| v_n\|_{\mu}^{q-1}\|v_n^+\|_{\frac{1}{2},2},\label{e:4.11}\\
    \int_M |u_n|^\frac{(q-1)p}{q}|v_n^+|dx&\leq& C\| u_n\|_{\mu}^\frac{(q-1)p}{q}\|v_n^+\|_{\frac{1}{2},2}.\label{e:4.12}
 \end{eqnarray}
 Moreover, it also holds that
 \begin{equation}\label{e:4.13}
   \int_M|u_n^+|dx\leq C\|u_n^+\|_{\frac{1}{2},2},\quad   \int_M|v_n^+|dx\leq C\|v_n^+\|_{\frac{1}{2},2}.
 \end{equation}
By (\ref{e:4.8})-(\ref{e:4.13}), we deduce
\begin{eqnarray}\label{e:4.14}
  \|z_n^+\|_E^2 &\le&   C(\| u_n\|_{\mu}^{p-1}+\| v_n\|_{\mu}^\frac{(p-1)q}{p}+1)\|u_n^+\|_{\frac{1}{2},2}+\nonumber \\
   && + C(\| v_n\|_{\mu}^{q-1}+\| u_n\|_{\mu}^\frac{(q-1)p}{q}+1)\|v_n^+\|_{\frac{1}{2},2}+\|z_n^+\|_E.
\end{eqnarray}
Hence this and (\ref{e:4.5}) lead to
\begin{eqnarray}\label{e:4.15}
   \|z_n^+\|_E^2 &\leq& C(\| z_n\|_E^\frac{p-1}{\mu}+\| z_n\|_E^\frac{(p-1)q}{\mu p}+1)\|z_n\|_E+\nonumber \\
   &+& C(\| z_n\|_E^\frac{q-1}{\mu}+\|z_n\|_E^\frac{(q-1)p}{\mu q}+1)\|z_n\|_E + \|z_n\|_E.
\end{eqnarray}
 For any $z^-\in E_-$,
then the similar arguments will lead to
\begin{eqnarray}\label{e:4.16}
  \|z_n^-\|_E^2 &\leq& C(\| z_n\|_E^\frac{p-1}{\mu}+\| z_n\|_E^\frac{(p-1)q}{\mu p}+1)\|z_n\|_E+\nonumber \\
  &+& C (\| z_n\|_E^\frac{q-1}{\mu}+\| z_n\|_E^\frac{(q-1)p}{\mu q}+1)\|z_n\|_E + \|z_n\|_E.
\end{eqnarray}
Adding (\ref{e:4.15}) and (\ref{e:4.16})   yields
\begin{eqnarray}\label{e:4.17}
 && \|z_n\|_E^2 \leq C(\| z_n\|_E^\frac{p-1}{\mu}+\| z_n\|_E^\frac{(p-1)q}{\mu p}+1)\|z_n\|_E+\nonumber \\
  &+& C(\| z_n\|_E^\frac{q-1}{\mu}+\| z_n\|_E^\frac{(q-1)p}{\mu q}+1)\|z_n\|_E
  +2 \|z_n\|_E.
\end{eqnarray}
By the assumptions  on $p,q,\mu$ above (${\bf H}_1$), it is easily checked that
 the total exponent each term in the right-hand side of (\ref{e:4.17}) is less than $2$. It follows
  that the sequence $\{z_n\}$ is bounded in $E$. Claim  is proved.
\par
Passing to a subsequence we may assume that for some $z\in E$,$z_n\rightharpoonup z$ weakly in  $E$.
From here on a usual argument based on the compactness of $D\mathcal{H}$ and invertibility of $B$ give the
existence of a subsequence of $z_n= B^{-1}(D\mathfrak{L}(z_n)+D\mathcal{H}(z_n))$ that converges
$z$ in $E$. So the $(PS)_c$-condition is verified.

\section{Proof of the Theorems}
The proof of Theorem \ref{the:1.1} is based on an application of the following theorem of Benci and Rabinowtitz\cite{BeRa}.
\begin{theorem}\label{the:5.1}
(Indefinite Functional Theorem). Let $H$ be a real Hilbert space with $H=H_1\bigoplus H_2$. satisfies the Palais-Smale conditon, and
\begin{flalign*}
(I_1):& \mathfrak{L}(z)=\frac{1}{2}(Lz,z)-\mathcal{H}(z), \text{where}\ L:H\rightarrow H \ \text{is bounded and self-adjoint},&\\
&\text{ and}\  L \ \text{leaves} \ H_1\  \text{and} \ H_2\  \text{invariant};&\\
(I_2):& D\mathcal{H} \ is \ compact;&\\
(I_3):& \text{there exists a subspace}\ \widetilde{H}\subset H\ \text{and sets}\ S\subset H, Q\subset \widetilde{H}\ \text{and constants}\ \alpha>\omega \ \text{such that} &\\
&(i): S\subset H_1 \ and \ \mathfrak{L}\mid_S \geq\alpha,&\\
&(ii): Q \ is \ bounded \ and\ \mathfrak{L}\leq \omega \ on \ the \ boundary \ \partial Q \ of \ Q \in \ \widetilde{H};&\\
&(iii): S \ and \ \partial Q \ link.&
\end{flalign*}
then $\mathfrak{L}$ possesses a critical value $\ c\geq \alpha$.
\end{theorem}
\par
Before giving the geometric conditions for the first linking property, We set $s_1$, $s_2$, $\rho>0$ with $0<\rho<s_2$, $[0,s_1 e^+]=:\{se^+|0\leq s\leq s_1\}$, and
$$
Q=[0,s_1 e^+]\oplus (\overline{B}_{s_2}\cap E_-),\quad \widetilde{H}=\text{span}[e^+]\oplus E_-,\quad S=\partial B_\rho \cap E_+,
$$
where $B_\rho$ denotes an open ball with radius $\rho$ centered at the origin, $e^+=(\xi^+,\eta^+)\in E_+$ with $\eta^+$ some eigenspinor of $P$ corresponding to the first positive eigenvalue $\lambda_1^{\text{CHI}}$.
\begin{lemma}\label{lem:5.2}
There exists  $\rho>0$ and $\alpha>0$ such that
$$
\mathfrak{L}(z)\geq\alpha\quad  \forall z\in S.
$$
\end{lemma}
\text{proof:} Conditon $(H_1)$,$(H_3)$ and $(H_4)$ imply that for any $\varepsilon>0$ there exists a constant $C(\varepsilon)>0$ such that

\begin{equation}\label{e:5.1}
H(x,u,v)\leq \varepsilon(\mid u\mid^2+\mid v\mid^2)+C(\varepsilon)(\mid u\mid^{p}+\mid v\mid^{q})
\end{equation}
for all  $(u,v)\in E$, combining (\ref{e:5.1}) and the Sobolev embedding, it is straightforward to show that
$$
\mathfrak{L}(z^+)\geq (\frac{1}{2}-C_4\varepsilon)\parallel z^+\parallel_E^2-C_5C(\varepsilon)(\parallel z^+\parallel_E^{p}+\parallel z^+\parallel_E^{q})
$$
for some constants $C_4>0$ and $C_5>0$. Thus we can fix $\varepsilon<\frac{1}{2C_4}$ and take enough small $\rho>0,\alpha>0$ such that $\mathfrak{L}(z^+)\geq\alpha$
on $S$.
\begin{lemma}\label{lem:5.3}
There exists  $s_1$, $s_2$, $\rho>0$ with $0<\rho<s_2$ such that
$$
\mathfrak{L}(z)\leq 0,\quad \forall z\in \partial Q.
$$
\end{lemma}
\text{proof:} Note that the boundary $ \partial Q$ of the cylinder $Q$ is taken in the space $\widetilde{H}$, and consists of three parts,
namely the bottom $Q\cap\{s=0\}$, the lid $Q\cap \{s=s_1\}$, and $[0,s_1e^+]\oplus(\partial B_{s_2}\cap E_-)$. Clearly $\mathfrak{L}(z)\leq 0$ on the bottom because by $(H_3)$. For the remaining two parts of the boundary we first observe that, for $z=z^-+re^+\in \widetilde{H}$,
\begin{equation}\label{e:5.2}
 \mathfrak{L}(z^-+re^+)=\frac{1}{2}r^2-\frac{1}{2}\parallel z^-\parallel_E^2- \mathcal{H}(z^-+re^+).
\end{equation}
 By definition of $E_{+}$ we have $\xi^+=\mid P\mid^{-1}P\eta^+=\eta^+$, therefore, $e^+=(\eta^+,\eta^+)$.

We set $z^-=(u^-,v^-)$, for $z^-+re^+=(u^-+r\eta^+,v^-+r\eta^+)$, using $(\ref{e:4.3})$, we have
\begin{equation}\label{e:5.3}
 \mathfrak{L}(z^-+re^+)\leq\frac{1}{2}r^2-\frac{1}{2}\parallel z^-\parallel^2_E- C_2\int_M\mid u^-+r\eta^+\mid^\mu-C_2\int_M\mid v^-+r\eta^+\mid^\mu+C_3.
\end{equation}
Thus, writing $v^-=t\eta^++\widehat{v}$, where $\eta^+$ is orthogonal to $\widehat{v}$ in $L^2(M,\Sigma M)$. By definition of $E_{\pm}$ we have
$$
 u^-=-\mid P\mid^{-1}Pv^-=-\mid P\mid^{-1}P(t\eta^++\widehat{v})=-t\eta^+-\mid P \mid^{-1}P\widehat{v},
$$
Similarly, $\eta^+$ is orthogonal to $\mid P\mid^{-1}P\widehat{v}$ in $L^2(M,\Sigma M)$. By the H$\ddot{o}$lder's inequality,
$$
(r+t)\int\mid \eta^+\mid^2dx=\int\langle r\eta^++v^-,\eta^+\rangle dx\leq\parallel r\eta^++v^-\parallel_{\mu}\parallel \eta^+\parallel_{\frac{\mu}{\mu-1}}
$$
which implies
\begin{equation}\label{e:5.4}
  r+t\leq C_6\parallel r\eta^++v^-\parallel_{\mu}
\end{equation}

for some constant $C_6$ depending on $\eta^+$. Similarly, we have
\begin{equation}\label{e:5.5}
  r-t\leq C_7\parallel r\eta^++u^-\parallel_{\mu}.
\end{equation}

Therefore, we deduce from (\ref{e:5.3}), (\ref{e:5.4}) and (\ref{e:5.5})that
\begin{equation}\label{e:5.6}
  \mathfrak{L}(z^-+re^+)\leq \frac{1}{2}r^2-Cr^\mu+C.
\end{equation}
By $\mu>2$, taking $r=s_1$ large enough we see in (\ref{e:5.6}) that $ \mathfrak{L}(z^-+re^+)<0$ on the lid $Q\cap \{s=s_1\}$.
\par
For $z^-+re^+ \in [0,s_1e^+]\oplus(\partial B_{s_2}\cap E_-)$, we deduce from  the condition $(H_3)$ that
\begin{equation*}
  \mathfrak{L}(z^-+re^+)\leq \frac{1}{2}r^2-\frac{1}{2}\parallel z^-\parallel_{E}^2.
\end{equation*}
 Taking $\parallel z^-\parallel_E=s_2$ large enough, it holds that
\begin{equation*}
  \mathfrak{L}(z^-+re^+)\leq 0.
\end{equation*}
The desired result is proved.
\par
\textbf{Proof of Theorem \ref{the:1.1}.}  Let $H=E, H_1=E_+, H_2=E_-$, we apply Lemma(\ref{lem:4.1}) to the functional $\mathfrak{L}$. The Palais-Smale conditon is satisfied. We can use the standard methods to show that Condition $I_1, I_2$ and $\mathfrak{L}$  is continuously differentiable. The geometric conditions $I_3(i),(ii)$ is proved in Lemma(\ref{lem:5.2}) and Lemma(\ref{lem:5.3}). For the proof of $(I_3)(iii)$ we refer to \cite{BeRa}. Therefore $\mathfrak{L}$ possesses a critical value point $z\in E$ and satisfies $\mathfrak{L}(z)\geq \alpha>0$.
\par
To obtain the theorem \ref{the:1.2}, we recall the Generalized fountain theorem for semi-definite functionals, see \cite{Fel} for the detailed exposition.
\begin{theorem}\label{the:5.4}
(Generalized fountain theorem) Let $X$ be a real Hilbert space with $X=Y\oplus Z$, where $Y$ is closed and seperable,  $Z= \overline{\bigoplus\limits_{j=0}^{\infty}\mathbb{R}\epsilon_j}$. $\mathfrak{L}\in C^1(E)$ is an even functional, i.e. $\mathfrak{L}(z)=\mathfrak{L}(-z)$ for all $z\in X$.
Assume $\mathfrak{L}$ is $\tau-$upper semicontinuous and $D\mathfrak{L}$ is weakly sequentially continuous. If for every $k\in \mathbb{N}$, there
exists $\rho_k>r_k>0$ such that
\begin{flalign*}
&(A_1): a_k:=\inf_{\substack{z\in Z_k\\ \parallel z\parallel= r_k}}\mathfrak{L}(z)\rightarrow \infty, k\rightarrow \infty; &\\
&(A_2): b_k:=\sup_{\substack{z\in Y_k\\ \parallel z\parallel=\rho_k}}\mathfrak{L}(z)\leq 0\  and\  d_k:=\sup_{\substack{z\in Y_k\\ \parallel z\parallel\leq\rho_k}}\mathfrak{L}(z)<\infty;&\\
&(A_3): \mathfrak{L} \ \text{satisfies the Palais-Smale condition};&
\end{flalign*}
where
\begin{equation*}
  Y_k:= Y\oplus\left( \bigoplus\limits_{j=1}^{k}\mathbb{R}\epsilon_j\right),\quad Z_k:= \overline{\bigoplus\limits_{j=k}^{\infty}\mathbb{R}\epsilon_j},
\end{equation*}
$\{\epsilon_j\}_{j=1}^\infty$ is a total orthonormal sequence in $Z$.
then $\mathfrak{L}$ has an unbounded sequence of critical values.
\end{theorem}
\par
We will define the following subsets for giving the geometric conditions of the linking property:
\begin{equation*}
  Y_k:= E_{-}\oplus\left( \bigoplus\limits_{j=1}^{k}\mathbb{R}e_j\right),\quad Z_k:= \overline{\bigoplus\limits_{j=k}^{\infty}\mathbb{R}e_j},
\end{equation*}
\begin{equation*}
   B_k:=\{z\in Y_k\mid \parallel z\parallel_E\leq \rho_k\}, N_k:=\{z\in Z_k\mid \parallel z\parallel_E=r_k\}, \partial B_k:=\{z\in Y_k\mid \parallel z\parallel_E= \rho_k\},
\end{equation*}
where $0<r_k<\rho_k,k\geq 2.$

\begin{lemma}\label{lem:5.5}
There exists  $\rho_k> r_k>0$ such that
\begin{flalign*}
(A_1)& a_k:=\inf_{\substack{z\in N_k}}\mathfrak{L}(z)\rightarrow \infty, k\rightarrow \infty;& \\
(A_2)& b_k:=\sup_{\substack{z\in \partial B_k}}\mathfrak{L}(z)\leq 0\  and\  d_k:=\sup_{\substack{z\in B_k}}\mathfrak{L}(z)<\infty.&
\end{flalign*}
\end{lemma}
Proof: (i) Let $z=(u,v)\in Z_k$, $T=\max\{p,q\}, t=\min\{p,q\}$, Then by (\ref{e:5.1}), which implies that
\begin{eqnarray}\label{e:5.7}
 \mathfrak{L}(z)&=&\frac{1}{2}\parallel z\parallel_E^2- \int_M H(x,z)dx\nonumber\\
  &\geq&\ \frac{1}{2}\parallel z\parallel_E^2-\varepsilon(\parallel u\parallel_2^2+\parallel v\parallel_2^2)-C(\varepsilon)(\parallel u\parallel_{p}^p+\parallel v\parallel_q^{q})\nonumber\\
 &\geq&\ (\frac{1}{2}-C_8\varepsilon)\parallel z\parallel_E^2-C_9C(\varepsilon)(\parallel z\parallel_E^{p}+\parallel z\parallel_E^{q})
\end{eqnarray}
where $C_8,C_9$ is constants.  Choosing $ \varepsilon=\frac{1}{4C_8}$, for $0< \alpha_k\leq 1$ and $\parallel z\parallel_E\geq1$, we obtain from (\ref{e:5.7}) that
\begin{eqnarray}\label{e:5.8}
 \mathfrak{L}(z)&\geq&\frac{1}{4}\parallel z\parallel_E^2-C(\varepsilon)C_9C_{10} \alpha_k^{p}\parallel z\parallel_E^{p}-C(\varepsilon)C_9C_{10}\alpha_k^{q}\parallel z\parallel_E^{q}\nonumber\\
  &\geq&\frac{1}{4}\parallel z\parallel_E^2-2C(\varepsilon)C_9C_{10} \alpha_k^{t}\parallel z\parallel_E^{T}
\end{eqnarray}
Then we have for $\parallel z\parallel_E=r_k:=(4TC(\varepsilon)C_9C_{10}\alpha_k^{t})^\frac{1}{2-T}$,
$$  \mathfrak{L}(z)\geq (\frac{1}{4}-\frac{1}{2T})(4TC(\varepsilon)C_9C_{10}\alpha_k^{t})^\frac{2}{2-T}$$
We know by lemma 3.8 in \cite{Wil} that $\alpha_k\rightarrow 0$ as $k\rightarrow \infty$, so $ \mathfrak{L}(z)\rightarrow \infty$ as $k\rightarrow \infty$,
and condition $(A_1)$ is satisfied.
\par
(ii) Let $z=z^-+w^+\in Y_k$, where $z^-=(z_1^-,z_2^-)\in E_-$, $w^+=(w_1^+,w_2^+)\in \bigoplus\limits_{j=1}^{k}\mathbb{R}e_j$.
The condition $(H_2)$ and $(H_4)$ implies that for every $\delta>0$ there exists a constant $C(\delta)>0$ such that
$$
H(x,z_1^-+w_1^+,z_2^-+w_2^+)\geq \delta(\mid z_1^-+w_1^+\mid^2+\mid z_2^-+w_2^+\mid^2)-C(\delta).
$$
Since $E_+$ is orthogonal to  $E_-$ in $L^2(M,\Sigma M\oplus \Sigma M)$, we deduce
\begin{eqnarray}\label{e:5.9}
 &-& \int_M H(x,z_1^-+w_1^+,z_2^-+w_2^+)dx\leq-\delta(\parallel z_1^-+w_1^+\parallel_2^2+\parallel z_2^-+w_2^+\parallel_2^2)+C(\delta)\mid M\mid\nonumber\\
  &\leq&-\delta(\parallel z_1^-\parallel_2^2+\parallel w_1^+\parallel_2^2+\parallel z_2^-\parallel_2^2+\parallel w_2^+\parallel_2^2)+C(\delta)\mid M\mid\nonumber\\
   &\leq&-\delta(\parallel w_1^+\parallel_2^2+\parallel w_2^+\parallel_2^2)+C(\delta)\mid M\mid
\end{eqnarray}
 All norms are equavalent on the space $\bigoplus\limits_{j=1}^{k}\mathbb{R}e_j$, there exists a constant $C>0$ such that
 $$
 C\parallel  w^+\parallel_E^2\leq\parallel  w^+\parallel_{L^2\times L^2}^2=\parallel w_1^+\parallel_2^2+\parallel w_2^+\parallel_2^2.
$$
Combining (\ref{e:3.1}) with (\ref{e:5.9}) yields
 \begin{eqnarray}\label{e:5.10}
 \mathfrak{L}(z)&=&-\frac{1}{2}\parallel z^-\parallel_E^2+\frac{1}{2}\parallel w^+\parallel_E^2- \int_M H(x,z)dx\nonumber\\
  &\leq&-\frac{1}{2}\parallel z^-\parallel_E^2+\frac{1}{2}\parallel w^+\parallel_E^2-\delta(\parallel w_1^+\parallel_2^2+\parallel w_2^+\parallel_2^2)+C(\delta)\mid M\mid\nonumber\\
  &\leq&-\frac{1}{2}\parallel z^-\parallel_E^2+\frac{1}{2}\parallel w^+\parallel_E^2-C(\delta)\parallel  w^+\parallel_E^2+C(\delta)\mid M\mid.
\end{eqnarray}

Takeing $\delta>\frac{1}{2C}$, this shows that $\mathfrak{L}(z)\rightarrow -\infty$ as $\parallel z\parallel_E\rightarrow \infty$, so  $(A_2)$ is satisfied for $\rho_k$
large enough. The desired result is proved.
\par
Let $\Pi_-:E\rightarrow E_-$ and $\Pi_+:E\rightarrow E_+$ be the orthogonal projections, $\{\theta_j\}_{j=1}^\infty$ be an orthonormal basis of $E_-$.
On $E$ we consider a new norm
\begin{equation*}
 \mid\mid\mid\mid z\mid\mid\mid\mid:=\max\left(\sum\limits_{j=0}^{\infty}\frac{1}{2^{j+1}}\mid (\Pi_-z,\theta_j)\mid,\parallel \Pi_+z\parallel_E\right),
\end{equation*}
We use the $\tau-$topology is generated by the norm $\mid\mid\mid\mid \cdot\mid\mid\mid\mid$ (see\cite{KrSz} or \cite{Wil}). It is clear that $\parallel \Pi_+z\parallel_E\leq\mid\mid\mid\mid z\mid\mid\mid\mid\leq\parallel z\parallel_E$. Moreover, if $\{z_n\}$ is a bounded sequence in $E$ then
$$z_n\xrightarrow{\tau} z\Leftrightarrow \Pi_-z_n\rightharpoonup \Pi_-z \ and \ \Pi_+z_n\rightarrow \Pi_+z.$$
\par
Let $\mathfrak{L}\in C^1(E)$, recall that $\mathfrak{L}$ is $\tau-$upper semicontinous if $z_n\xrightarrow{\tau} z$ implies $\mathfrak{L}(z)\geq\overline{\lim\limits_{n\to \infty}}\mathfrak{L}(z_n)$.

\begin{lemma}\label{lem:5.6}
Under assumption $(H_1)$ and $(H_3)$, $\mathfrak{L}$ is $\tau-$ upper semicontinuous and $D\mathfrak{L}$ is weakly sequentially continuous.
\end{lemma}
Proof:(i) Assume $z_n\in E$ such that $z_n\xrightarrow{\tau} z$ and $\mathfrak{L}(z_n)\geq c$. By the definition of $\tau$
 we obtain $ \Pi_+z_n\rightarrow \Pi_+z$. Noting that $\mathfrak{L}(z_n)\geq c$ and $H\geq 0$, then $\Pi_-z_n$ is bounded and
$ \Pi_-z_n\rightharpoonup \Pi_-z$. Now since the embedding $E\hookrightarrow L^p\times L^q$ is compact, then $z_n\rightarrow z$ in  $L^p\times L^q$. Thus
up to a subsequence $z_n\rightarrow z$ a.e. on $M$, and by $H\in C^0(\Sigma M\oplus \Sigma M, \mathbb{R})$ is $C^1$ in the fiber direction, we obtain $H(x,z_n)\rightarrow H(x,z)$ a.e. on $M$. Using the Fatou lemma
and the weak lower semicontinuity of the norm $\parallel\cdot \parallel_E$ that $\mathfrak{L}(z)\geq \overline{\lim\limits_{n\to \infty}}\mathfrak{L}(z_n)\geq c$. Therefore, $\mathfrak{L}$ is $\tau-$upper semicontinuous.
\par
(ii) Assumen $(u_n,v_n)\in E$ such that $(u_n,v_n)\rightharpoonup (u,v)$, then by Rellich's theorem $u_n\rightarrow u$ in  $L^p(M,\Sigma M)$ and $v_n\rightarrow v$ in  $L^q(M,\Sigma M)$.   then by the H\"{o}lder inequality we have

\begin{eqnarray*}
 && \mid D\mathfrak{L}(u_n,v_n)(\xi,\zeta)-D\mathfrak{L}(u,v)(\xi,\zeta)\mid\leq \\
  &\leq&\int_M \mid \xi\mid\mid H_u(x,u_n,v_n)-H_u(x,u,v)\mid dx+\int_M \mid \zeta\mid\mid H_v(x,u_n,v_n)-H_v(x,u,v)\mid dx\\
   &\leq&\parallel\xi\parallel_p\parallel H_u(x,u_n,v_n)-H_u(x,u,v)\parallel_\frac{p}{p-1}+\parallel\zeta\parallel_q\parallel H_v(x,u_n,v_n)-H_v(x,u,v)\parallel_\frac{q}{q-1}
\end{eqnarray*}
It follows from
$$L^p(M,\Sigma M)\times L^q(M,\Sigma M)\rightarrow L^\frac{p}{p-1}(M,\Sigma M)$$
 and
 $$L^p(M,\Sigma M)\times L^q(M,\Sigma M)\rightarrow L^\frac{q}{q-1}(M,\Sigma M)$$
 are continuous  by $(H_1)$ that
$$
 D\mathfrak{L}(u_n,v_n)(\xi,\zeta)\rightarrow D\mathfrak{L}(u,v)(\xi,\zeta)\qquad \text{as}\ n\rightarrow \infty.
$$
This shows that $D\mathfrak{L}$ is weakly sequentially continuous.
\par

Now we use Theorem \ref{the:5.4} to obtain infinitely many critical points of the even functional $\mathfrak{L}$ in Theorem \ref{the:1.2}.
\par
\textbf{Proof of Theorem \ref{the:1.2}}  We know by Lemma \ref{lem:4.1} that $\mathfrak{L} $ satisfies the Palais-Smale condition, Lemma \ref{lem:5.5} gives the geometric
conditions $A_1$ and $A_2$ of the linking geometry. Lemma \ref{lem:5.6} implies $\mathfrak{L}$ is $\tau-$ upper semicontinuous and $D\mathfrak{L}$ is weakly sequentially continuous. Then  $\mathfrak{L} $ has an unbounded sequence of critical values.

$\bf{Founding}$ The *corresponding author: Xu Yang was supported by the NSFC
(grant no.11801499) of China.

%
%


\begin{thebibliography}{5}
%
\bibitem{Abb} A. Abbondandolo, A new cohomology for the Morse theory of strongly indefinite functionals on Hilbert spaces.
{\it Topol. Methods NOnlinear Anal.},{\bf 9}(1997), no.2,325-382.
\bibitem{Ada} R. Adams, Sobolev Space.{\it Academic Press, New York}(1975).

\bibitem{Am1} B. Ammann, A variational Problem in Conformal Spin Geometry.{\it Journal of Geometry  Physics} ,{\bf62}(2003),no.2 , :213-223

\bibitem{Am2} B. Ammann, The smallest Dirac eigenvalue in a spin-conformal class and cmc-immersions.
{\it Comm. Anal. Geom.}, {\bf 17}(2009), 429-479.

\bibitem{AnVo} S. Angenent and R.van der Vorst, A superquadratic indefinite elliptic system and its Morse-Conley-Floer homology.
{\it Math. Z.},{\bf 231}(1999), no.2, 203-248.

\bibitem{BaCo} C. J. Batkam, F. Colin,Generalized fountain theorem and applications to strongly indefinite semilinear problems,
{\it J.Math. Anal. appl.}, {\bf 405}(2013), 438-452.
\bibitem{BaDi1} T. Bartsch and Y. Ding, Homoclinic solutions of an infinite-dimensional Hamiltonian system.
{\it Math. Z.}, {\bf 240}(2002), 289-310.

\bibitem{BaDi2} T. Bartsch and Y. Ding, Periodic solutions of superlinear beam and membrane equations with perturbations from symmetry.
{\it Nonlinear Analyss}, {\bf 44}(2001),727-748.

\bibitem {BaLi} A. Bahri and P. L. Lions, Solutions of superlinear elliptic equations and their Morse indices.
{\it Commun.Pure Appl.Math.}, {\bf 45}(1992),1205-1215.

\bibitem{BeRa} V. Benci and P. H. Rabinowitz, Periodic solutions of Hamiltonian systems.
{\it Pure Appl. Math.}, {\bf 31}(1978),157-184.

\bibitem{CJLW1} Q. Chen, J. Jost,J. Li and G.Wang,Dirac-harmonic maps,{\it Math. Z.},{\bf 254}(2006),409-432.

\bibitem{CJLW2} Q. Chen, J. Jost,J. Li and G.Wang,Nonlinear Dirac equations on Riemann surfaces,{\it Ann. Global Anal.Geom.},{\bf33}(2008),253-270.
\bibitem{DiL} Y. H. Ding, J. Y. Li,A boundary value problem for the nonlinear Dirac equation on compact spin manifold,{\it Calc. Var.},{\bf 57}(2018).
\bibitem{Fel} P. Felmer, Periodic solutions of 'supequadratic' Hamiltonian systems, {\it J. Differential Euations},{\bf 102}(1993),188-207.

\bibitem{FaSc} S. Farinell, G. Schwarz, On the spectrum of the Dirac operator under boundary conditions, {\it J. Geom. Phys.},{\bf 28}(1998),67-84.

\bibitem{FF} D.G. De Figueiredo, P.L. Felmer,  On superquadratic elliptic systems. {\it Trans. Am.
Math. Soc.}, {\bf 343}(1994), 97–116.
\bibitem{Fri1} T. Friedrich, On the spinor representation of surfaces in Eulidean 3-space, {\it J. Geom. Phy.},{\bf 28}(1998),143-157.
\bibitem{Fri2} T. Friedrich, \emph{Dirac Operators in Riemannian Geometry}, Grad. Stud. Math., vol. 25, Amer. Math. Soc.,Providence, RI, 2000.
\bibitem{GoL1} W. Gong and G. Lu, On Dirac equation with a potential and critical Sobolev exponent.
{\it Commun. Pure Appl. Anal.} {\bf 14} (2015), 2231 - 2263.
\bibitem{GoL2} W. Gong and G. Lu, Existence results for coupled Dirac systems via Rabinowitz-Floer theory,
DOI:10.1007/s11784-013-0116-5.
\bibitem{HMR} O. Hijazi, S. Montiel,A. Rold$\acute{a}$n,  Eigenvalue boundary problems for the Dirac operator. {\it Commun. Math.Phys.},{\bf 231}(2002), 375–390.
\bibitem{HMZ} O. Hijazi, S. Montiel,X. Zhang,$\acute{a}$n,  Eigenvaluevs of the Dirac operator on manifolds with boundary. {\it Commun. Math.Phys.},{\bf 221}(2001), 255–265.
\bibitem{HV} J. Hulshof, R. van der Vorst, Differential systems with strongly indefinite variational
structure. {\it J. Funct. Anal.}, {\bf 114}(1993), 32–58.

\bibitem{Iso1} T. Isobe, Existence results for solutions to nonlinear Dirac equations on compact spin manifolds. {\it Manuscripta math}, {\bf 135} (2011), 329 - 360.
\bibitem{Iso2} T. Isobe, Nonlinear Dirac equations with critical nonlinearities on compact spin manifolds. {\it J.Funct. Anal.}, {\bf 260} (2011), 253 - 307.
\bibitem{Iso3} T. Isobe, A perturbation method for spinorial Yamabe type equations on $S^m$ and its application. {\it Math. Ann.}, {\bf 355} (2013), 1255-1299.
\bibitem{KrSz} W. Kryszewski and A. Szulkin, An infinite-dimensional Morse theory with applications. {\it Trans. Amer. Math. Soc.}, {\bf 349}(1997), 3181-3234.
\bibitem{LaM} H.B. Lawson and M.L. Michelson, \emph{Spin Geometry}. Princeton University Press, 1989.

\bibitem{Maa} A. Maalaoui, Rabinowitz-Floer homology for superquadratic Dirac equations on spin manifolds. {\it  J. Fixed Point Theory Appl.}, {\bf 13} (2013), 175-199.

\bibitem{MaM} A. Maalaoui and V. Martino, The Rabinowitz-Floer homology for a class of semilinear problems and applications. {\it J. Funct. Anal.}, {\bf 269} (2015), 4006-4037.
 \bibitem{Rau} S. Raulot, A Sobolev-like inequality for the Dirac operator, {\it J. Funct. Anal.},{\bf 256}(2009),1588-1617.

\bibitem{Str} M. Struwe, Vatiational Methods.  Applications to Nonlinear Partial Differential Equations and Hamiltonian Systems, 4th edn. {\it Springer,Berlin}(2008).

\bibitem{Wil} M. Willem, Minimax theorems.  Progress in Nonlinear Partial Differential Equations and Their Applications.{\it  Birkh$\ddot{a}$user, Boston } (1996).

\bibitem{YJL} X. Yang, R. Jin, G. Lu, Solutions of Dirac equations on compact manifolds
via saddle point reduction. {\it  J. Fixed Point Theory Appl.} ,{\bf 19}(2017),215-229.

\bibitem{Yan} X. Yang, Existence results for solutions to nonlinear Dirac systems on compact spin manifolds. {\it  Advanced Nonlinear Studies} ,DOI: 10.1515/ans-2017-6034.

\end{thebibliography}
\end{document}